\documentclass[12pt,oneside]{amsart}
\usepackage{amsaddr}
\usepackage[utf8]{inputenc}
\usepackage[T1]{fontenc}
\usepackage{lmodern}
\usepackage{amssymb}
\usepackage{geometry}
\usepackage{amsmath}
\usepackage{amsthm}
\usepackage{graphicx}
\usepackage{color}
\usepackage[normalem]{ulem}
\usepackage[onehalfspacing]{setspace}
\usepackage{verbatim}
\usepackage[hidelinks]{hyperref}

\newcommand{\IR}{\ensuremath{\mathbb{R}}}
\newcommand{\IN}{\ensuremath{\mathbb{N}}}

\newcommand{\IC}{\ensuremath{\mathbb{C}}}

\newcommand{\IP}{\ensuremath{\mathbb{P}}}
\newcommand{\IE}{\ensuremath{\mathbb{E}}}

\newcommand{\norm}[1]{\left\Vert#1\right\Vert}
\newcommand{\set}[1]{\left\{#1\right\}}

\newcommand{\brackets}[1]{\left(#1\right)}

\newcommand{\diag}{\mathop{\mathrm{diag}}}

\newtheorem{thm}{Theorem}

\theoremstyle{plain}

\newtheorem{prop}{Proposition}
\theoremstyle{definition}

\title[On the worst-case error of least squares algorithms]{
On the worst-case error of \\
least squares algorithms 
for $L_2$-approximation \\
with high probability}

\author{
Mario Ullrich
}

\address{Institut f\"ur Analysis, 
Johannes Kepler Universit\"at Linz, Austria}
\email{
mario.ullrich@jku.at}


\date{\today}

\begin{document}

\begin{abstract}
It was recently shown in \cite{KU19} that, for $L_2$-approximation 
of functions from a Hilbert space, function values are almost 
as powerful as arbitrary linear information, 
if the approximation numbers are square-summable. 
That is, we showed that
\[
 e_n \,\lesssim\, \sqrt{\frac{1}{k_n} \sum_{j\geq k_n} a_j^2}
\qquad \text{ with }\quad k_n \asymp \frac{n}{\ln(n)},
\]
where 
$e_n$ are the sampling numbers and 
$a_k$ are the approximation numbers. 
In particular, if $(a_k)\in\ell_2$, then $e_n$ and $a_n$ are of the 
same polynomial order. 
For this, we presented an explicit (weighted least squares) algorithm 
based on i.i.d.~random points and proved that this works 
with positive probability. 
This implies the existence of a good deterministic sampling algorithm.

Here, we present a modification of the proof in \cite{KU19} 
that shows that the same algorithm works 
with probability at least $1-{n^{-c}}$ for all $c>0$.
\end{abstract}

\maketitle


Let $H$ be a Hilbert space of 
real- or complex-valued functions on a set $D$
such that point evaluation
$$
 \delta_x\colon H \to \IR,\quad f\mapsto f(x)
$$
is a continuous functional for all $x\in D$, 
which are usually called 
\emph{reproducing kernel Hilbert spaces}.
We consider numerical approximation of functions from such spaces, 
using only function values.
We measure the error in the space $L_2=L_2(D,\mathcal{A},\mu)$
of square-integrable functions with respect to an arbitrary 
measure $\mu$ such that $H$ is embedded into $L_2$.
This means that 
$H$ consists of square-integrable functions
such that two functions 
that are equal \mbox{$\mu$-almost} everywhere
are also equal point-wise.

We are interested in the \emph{$n$-th minimal worst-case error}
\[
e_n \,:=\, e_n(H) \,:=\, 
\inf_{\substack{x_1,\dots,x_n\in D\\ \varphi_1,\dots,\varphi_n\in L_2}}\, 
\sup_{f\in H\colon \|f\|_H\le1}\, 
\Big\|f - \sum_{i=1}^n f(x_i)\, \varphi_i\Big\|_{L_2},
\]
which is the worst-case error of an optimal algorithm that  
uses at most $n$ function values.
These numbers are sometimes called \emph{sampling numbers}. 
We want to compare $e_n$ with the \emph{$n$-th approximation number} 
\[
a_n \,:=\, a_n(H) \,:=\,
\inf_{\substack{L_1,\dots,L_n\in H'\\ \varphi_1,\dots,\varphi_n\in L_2}}\, 
\sup_{f\in H\colon \|f\|_H\le1}\, 
\Big\|f - \sum_{i=1}^n L_i(f)\, \varphi_i\Big\|_{L_2},
\]
where $H'$ is the space of all bounded, linear functionals on $H$.
This is the worst-case error of an optimal algorithm that uses 
$n$ linear functionals as information, 
and it is known that it equals the $n$-th singular value of 
the embedding ${\rm id}\colon H\to\!L_2$.
%
For an exposition of such approximation problems we refer to 
\cite{NW08,NW10,NW12}, especially~\cite[Chapter~26~\&~29]{NW12}, 
and references therein.
The main result of~\cite{KU19} is stated as follows.

\begin{thm}[\cite{KU19}]\label{thm:main}
 There are absolute constants $C,c>0$ and a sequence of natural numbers
 $(k_n)$ with $k_n\ge c n/\ln(n+1)$ such that the following holds.
 For any~$n\in\IN$, any measure space $(D,\mathcal A,\mu)$ and 
 any reproducing kernel Hilbert space $H$ of real-valued functions on $D$
 that is embedded into $L_2(D,\mathcal A,\mu)$, we have
 \[
  e_n(H)^2 \,\le\, \frac{C}{k_n} \sum_{j\geq k_n} a_j(H)^2.
 \]
\end{thm}


We refer to~\cite{KU19} for a thematic classification, 
further literature, 
and the implications for some long-standing open problems in the field.\\
This theorem was extended in~\cite{KUV19} to (complex) 
Hilbert spaces that may not be embedded into $L_2$, 
which may happen if the support of $\mu$ is not equal to $D$. 
Moreover, we've learned form \cite{KUV19} about the more recent 
paper~\cite{O10}, which allows for bounds on the singular values of 
random matrices with explicit constants. 
Combining~\cite{O10} with the proof technique from~\cite{KU19} we 
do not only see explicit constants from~\cite{KUV19}. 
We obtain that the method described below works with high probability, 
i.e., with probability at least $1-{n^{-c}}$ for all $c>0$.

\smallskip
\goodbreak

Before we state the main result,  
let us recall the method from~\cite{KU19}:\\
First of all, 
let 
${\rm id}\colon H\to L_2$ 
be the (injective) embedding form $H$ to $L_2(D,\mathcal{A},\mu)$, 
and $W:={\rm id}^*{\rm id}$. 
Since $W$ is positive and compact, 
there is an orthogonal basis $\mathcal B=\set{b_k \colon k\in\IN}$ of $H$ 
that consists of eigenfunctions of $W$.
Without loss of generality, we may assume that $H$ is 
infinite-dimensional. 
It is easy to verify that $\mathcal B$ is also orthogonal in $L_2$.
We may assume that the eigenfunctions are normalized in $L_2$
and that $\|b_1\|_H \leq \|b_2\|_H \leq \dots$, such that 
$a_k(H)=\|b_{k+1}\|_H^{-1}$.

\smallskip

Now let 
$k\in\IN$ 
(to be specified later), 
$x_1,\dots,x_n\in D$ be some given sampling nodes, 
and $V_k:={\rm span}\{b_1,\dots,b_k\}$. 
We then consider the algorithm 
\begin{equation}\label{eq:alg}
A_{n,{k}}(f) \,:=\, \underset{g\in V_{k}}{\rm argmin}\, \sum_{i=1}^n \frac{\vert g(x_i) - f(x_i) \vert^2}{\varrho_k(x_i)}, 
\end{equation}
where $\rho=\rho_k$ 
is given by
\[
 \varrho_k\colon D\to \IR, \quad \varrho_k(x) = \frac12 \left(
 \frac1k \sum_{j< k} b_{j+1}(x)^2  +  \frac{1}{\sum_{j\geq k} a_j^2} \sum_{j\geq k} a_j^2 b_{j+1}(x)^2
 \right).
\]
Note that, under mild assumptions, 
we have $A_{n,k}(f)=f$ whenever $f\in V_k$.\\
The \emph{worst-case error} of $A_{n,k}$ is defined as
\[
 e(A_{n,k},H) \,:=\, \sup_{f\in H\colon \|f\|_H\le1}\,  \big\|f - A_{n,k}(f)\big\|_{L_2}, 
\]
and we have $e_n(H)\leq e(A_{n,k},H)$ for every choice of $k$ and $x_1,\dots,x_n$.

In~\cite{KU19} we proved that, if $x_1,\dots,x_n$ are 
i.i.d.~random points with $\mu$-density $\rho$, 
then $e(A_{n,k_n})$ 
with $k_n \asymp n/\log(n)$
satisfies the bound in Theorem~\ref{thm:main} 
with positive probability. 
Here, we show that this holds with probability tending to 1, 
and we determine some explicit constants. 
(We did not try to optimize them.)
Roughly speaking, this shows that, asymptotically, almost all 
point sets lead to an algorithm for $L_2$-approximation that 
satisfies the bound above. This may increase the belief 
in the conjecture $e_n\asymp a_n$, 
see e.g.~\cite[Open Problem~140]{NW12}


Our improved result reads as follows.

\begin{thm}\label{thm:new}
For $n\ge2$ and $c>0$, let 
\[
k_n \,:=\, 2\cdot \left\lfloor \frac{n}{2^8\,(2+c)\,\ln(n)} \right\rfloor.
\] 
Then, for any measure space $(D,\mathcal A,\mu)$ and 
any reproducing kernel Hilbert space~$H$ of real- or 
complex-valued functions on $D$
 that is embedded into $L_2(D,\mathcal A,\mu)$, 
we have
 \[
  e_n(A_n,H)^2 \,\le\, \frac{4}{k_n}\, \sum_{j\geq k_n/2} a_j(H)^2
 \]
with probability at least $1-\frac8{n^c}$, 
where $A_n=A_{n,k_n}$ from~\eqref{eq:alg}.
\end{thm}

\bigskip
\goodbreak

\section*{The Proof}

The proof of Theorem~\ref{thm:new} is
almost the same as given in~\cite{KU19}, 
and is therefore very much inspired by 
the general technique to 
assess the quality of \emph{random information}
as developed in~\cite{HKNPU19b,HKNPU19a}. 
See also the references collected there.
In fact, we only replace 
\cite[Proposition~1]{KU19} (which is~\cite[Thm.~2.1]{MP06}) 
by~\cite[Lemma 1]{O10} 
to bound the singular values of the random matrices 
under consideration. \\
Let us note that the proof looks rather elementary, 
and it might be surprising that the results 
presented in~\cite{KU19} (and here), are not known for some time.
However, the way of controlling the 
'infinite-dimensional part' by adjusting the density $\rho$ accordingly, 
was seemingly invented in~\cite{KU19}, 
and this turned out to be essential.

First, let us note that the algorithm from~\eqref{eq:alg} 
can be written as
\[
A_{n,k}(f)=\sum_{j=1}^k (G^+ N f)_j b_j,
\]
where $N\colon H\to \IR^n$ with $N(f)=(\varrho(x_i)^{-1/2}f(x_i))_{i\leq n}$ 
is the weighted \emph{information mapping} and
$G^+\in \IR^{k\times n}$ is the Moore-Penrose inverse of the matrix 
\[
 G=(\varrho(x_i)^{-1/2} b_j(x_i))_{i\leq n, j\leq k} \in \IR^{n\times k},
\]
assuming that $G$ has full rank.

To give an upper bound on $e(A_{n,k})$,
let us assume that $G$ has full rank. 
For any $f\in H$ with $\Vert f\Vert_H\leq 1$, we 
let $P_k f$ be the orthogonal projection of $f$ to $V_k$, and obtain
\[\begin{split}
 \norm{f-A_{n,k}(f)}_{L_2}^2 \,&\le\, a_k^2 + \norm{P_k  f - A_{n,k}(f)}_{L_2}^2
 \,=\, a_k^2 + \norm{A_{n,k}(f- P_k f)}_{L_2}^2 \\
 &=\, a_k^2 + \norm{G^+ N(f- P_k f)}_{\ell_2^k}^2 \\
 &\le\, a_k^2 +\norm{G^+\colon \ell_2^n \to \ell_2^k}^2 \norm{N\colon P_k(H)^\perp \to \ell_2^n}^2.
\end{split}
\]
We've used $A_{n,k}(f)\in V_{k}$ in the first inequality, 
and $A_{n,k}(f)=f$ for $f\in V_k$ in the equality thereafter.
The norm of $G^+$ is the inverse of the $k$th largest (and therefore the smallest) singular value of the matrix $G$.
The norm of $N$ is the largest singular value of the matrix 
\[
 \Gamma =\big(\varrho(x_i)^{-1/2} a_j b_{j+1}(x_i)  \big)_{1\leq i \leq n, j\geq k} \in \IR^{n\times \infty}.
\]
To see this, note that 
$f\,=\,\sum_{j=1}^\infty \langle f,b_j\rangle_{L_2}\, b_j$ 
converges in $H$ for every $f\in H$, and therefore also point-wise, 
and that $\|f\|_H^2=\sum_{j=0}^\infty a_j^{-2}|\langle f,b_{j+1}\rangle_{L_2}|^2$.
Hence, $N=\Gamma \Delta$ on $P_k(H)^\perp$, 
where the mapping $\Delta\colon P_k(H)^\perp\mapsto\ell_2$ with 
$\Delta f=\left(\frac{\langle f,b_{j+1}\rangle_{L_2}}{a_j}\right)_{j\ge k}$
is an isomorphism.
This yields 
\vspace{-2mm}
\begin{equation}
\label{eq:basic}
 e(A_{n,k})^2 \,\leq\, a_k^2 + \frac{s_{\rm max}(\Gamma)^2}{s_{\rm min}(G)^2}.
\end{equation}

It remains to bound $s_{\rm min}(G)$ from below and $s_{\rm max}(\Gamma)$ from above.
Clearly, any nontrivial lower bound on $s_{\rm min}(G)$ automatically yields that the matrix $G$ has full rank.
To state our results, let
\[
\beta_k 
\,:=\, \brackets{\frac{1}{k} \sum_{j\geq k} a_j^2}^{1/2} 
\qquad\text{ and }\qquad
\beta'_k \,:=\, \beta_ {\lfloor k/2 \rfloor}.
\]
Note that $a_{2k}^2\le\frac1k(a_k^2+\hdots+a_{2k}^2)\le \beta_{k}^2$\, 
for all $k$ and thus $\max\{a_{k},\beta_k\} \leq \beta'_k$.

The rest of the paper is devoted to the proof of the following two claims: \\[2mm]
For each $k\le \frac{n}{128\cdot(2+c)\cdot\log(n)}$,
we have 
%

\begin{itemize}
\itemindent=15mm
\itemsep=2mm
	\item[{\bf Claim 1:}] \quad 
	$\displaystyle 
	\IP\Big(s_{\rm max}(\Gamma)^2 \,\leq\, n\,\frac{3(\beta'_k)^2}{2} \Big) \,\ge\, 1-\frac4{n^c}
	$
	\item[{\bf Claim 2:}] \quad\qquad 
$\displaystyle
\IP\Big(s_{\rm min}(G)^2 \,\geq\, \frac{n}{2} \Big) \,\;\ge\;\, 1-\frac4{n^c}
$
\end{itemize}
\medskip

Together with \eqref{eq:basic} and a union bound, this yields 
\[
 e(A_{n,k})^2 \,\le\, a_k^2 + 3(\beta'_k)^2 
 \,\leq\, 4 \beta_{\lfloor k/2 \rfloor}^2
\]
with probability at least $1-\frac8{n^c}$, 
which is the statement of Theorem~\ref{thm:new}.

\medskip
Both claims are based on~\cite[Lemma 1]{O10}, which we state here 
in a special case, i.e., we set $\delta=\frac4{n^c}$, see~\cite[top of p.~205]{O10}.
By $\norm{M}$ we denote the spectral norm of $M$.

\begin{prop}\label{prop:O}
Let $X$ be a random vector in $\IC^k$ or $\ell_2$ with  
$\|X\|_2\le R$ with probability 1, and let $X_1,X_2,\dots$ 
be independent copies of $X$. 
Additionally, let $E:=\IE(X X^*)$ satisfy $\|E\|\le 1$, 
and define 
\[
g(n,R,c) \,:=\, 4R\, \sqrt{\frac{(2+c)\ln(n)}{n}}.
\]
If $g(n,R,c)\le2$, then
\[
\IP\left(\bigg\|\sum_{i=1}^n X_i X_i^* - nE\bigg\|
	\,\le\, n\cdot g(n,R,c)\right)
\,\ge\, 1-\frac{4}{n^c}.
\]
\end{prop}
\medskip

\medskip

\begin{proof}[{\bf Proof of Claim 1}]
Consider independent copies $X_1,\hdots,X_n$ of the vector 
\[
X\,=\, \frac{1}{\beta'_k \sqrt{\varrho(x)}} 
	\Bigl(a_k b_{k+1}(x), a_{k+1} b_{k+2}(x), \hdots\Bigr)^\top,
\]
where $x$ is a random variable on $D$ with density $\varrho$.
Clearly, $\sum_{i=1}^n X_i X_i^* = \frac{1}{(\beta'_k)^2}\Gamma^* \Gamma$ 
with $\Gamma$ from above.
First observe
\[
 \norm{X}_2^2 \,=\, \frac{1}{(\beta'_k)^2 \varrho(x)} \sum_{j\geq k} a_j^2\, b_{j+1}(x)^2 
 \,\leq\, \frac{2}{(\beta'_k)^2} \sum_{j\geq k} a_j^2
 \,=\, 2k \,=:\, R^2.
\]
Since $E=\IE(X X^*)=\diag(\frac{a_{k}^2}{(\beta'_k)^2}, \frac{a_{k+1}^2}{(\beta'_k)^2}, 
	\hdots)$, 
we have $\|E\|=\frac{a_{k}^2}{(\beta'_k)^2}\le1$. 

Using $k\le \frac{n}{128(2+c) \ln(n)}$ and 
\[
g(n,R,c) \,=\, g\left(n, \sqrt{2k},c\right) 
\,=\, 
\sqrt{32\,(2+c)\, k \frac{\ln(n)}{n}} \le \frac12,
\] 
we obtain from 
Proposition~\ref{prop:O} that 
\[
 \IP\Big(\norm{\sum_{i=1}^n X_i X_i^* - nD} \le \frac{n}{2}\Big)
 \ge 1-\frac{4}{n^c}.
\]
This implies
\[\begin{split}
 s_{\rm max}(\Gamma)^2 &= \norm{\Gamma^*\Gamma} 
= (\beta'_k)^2 \norm{\sum_{i=1}^n X_i X_i^*} 
\leq (\beta'_k)^2\left(\norm{nE} + \norm{\sum_{i=1}^n X_i X_i^* - nE}\right) \\
&\leq n\, a_k^2 + n \frac{(\beta'_k)^2}{2}
\end{split}\]
with probability at least $1-4/n^c$ for all $k\le \frac{n}{128(2+c) \ln(n)}$.
This yields Claim~1.
\end{proof}

\begin{proof}[{\bf Proof of Claim 2}]
Consider $X=\varrho(x)^{-1/2}(b_1(x), \hdots, b_k(x))^\top$
with $x$ distributed according to $\varrho$.
Clearly, $\sum_{i=1}^n X_i X_i^* = G^*G$ with $G$ from above.
First observe
\[
 \norm{X}_2^2 \,=\, \varrho(x)^{-1} \sum_{j\le k} b_j(x)^2 
 \,\leq\, 2 k \,=:\, R^2.
\]
Since $E=\IE(X X^*)=\diag(1, \hdots,1)$ 
we have $\|E\|=1$. \\
Again, for 
$k\le \frac{n}{128(2+c) \ln(n)}$, we obtain 
$
g(n,R,c) \,\le\, \frac12,
$ 
and therefore, from 
Proposition~\ref{prop:O}, that 
$
 \IP\Big(\norm{\sum_{i=1}^n X_i X_i^* - nE} \le \frac{n}{2}\Big)
 \ge 1-\frac{4}{n^c}.
$
This implies that 
\[
 s_{\rm min}(G)^2 = s_{\rm min}(G^*G) \,\geq\, s_{\rm min}(nE) - \|G^*G - nE\|
 \,\geq\, n/2
\]
with probability at least $1-\frac4{n^c}$, 
and yields Claim~2.
\end{proof}

\subsection*{Acknowledgement}
I thank David Krieg for the  
indispensable discussions. \\
He decided not to be a co-author, because he thinks that he didn't 
contribute.

\bigskip


\end{document}